\RequirePackage{fix-cm}
\documentclass[smallextended]{svjour3}       
\smartqed  
\usepackage{graphicx}
\usepackage{amsmath,amsfonts,bm,mathtools}
\begin{document}

\title{Global solutions to general polynomial benchmark optimization problems
}


\author{Xiaojun Zhou         \and
        David Yang Gao       \and
        Chunhua Yang 
}


\institute{Xiaojun Zhou \at
              School of Science, Information Technology and Engineering, University of Ballarat, Victoria 3353, Australia.\\
              School of Information Science and Engineering, Central South University, Changsha 410083, China.\\
           \and
           David Yang Gao \at
              School of Science, Information Technology and Engineering, University of Ballarat, Victoria 3353, Australia.\\
           \and
           Chunhua Yang \at
              School of Information Science and Engineering, Central South University, Changsha 410083, China.
}

\date{Received: date / Accepted: date}

\maketitle

\begin{abstract}
The goal of this paper is to solve a class of high-order polynomial benchmark optimization problems, including the Goldstein-Price problem and the Three Hump Camel Back problem.
By using a generalized canonical duality theory, we are able to transform the nonconvex primal problems to concave dual problems over convex domain(without duality gap), which can be solved easily to obtain global solutions.
\end{abstract}

\section{Introduction}
\label{intro}
Polynomial optimization problems have been widely studied in various fields such as nonlinear algebra, semidefinite programming, and operations research, with extensive applications in production planning, location and distribution, engineering design, risk management, water treatment and distribution, chemical process design, pooling and blending, structural design, signal processing, robust stability analysis, design of chips, and much more (see \cite{Ref1,Ref2}).\\
\indent Due to the nonconvexity, traditional direct methods for solving polynomial optimization problems are usually very difficult, or even impossible. For example, the algebraic method for the task is to find all of the critical points firstly and then to identify the global minimizer(s) among all these critical points. This approach becomes inefficient when there exist numerous local minima. Also, linearization and relaxation techniques were used to compute an approximate optimal solution of the primal problem. However, the approximate optimal solution was not guaranteed to be the actual global optimum \cite{Ref3,Ref4}. In \cite{Ref5}, the so-called Z-eigenvalue methods were proposed to solve the best rank-one approximation problem, but they can be applied only for third-order polynomials. Besides these deterministic methods, stochastic techniques have also made significant contributions to the optimization applications of this kind \cite{Ref6,Ref7}. For example, the evolutionary computation method, could solve general problems in low dimension, but it failed to do well for large scale ones \cite{Ref8,Ref9,Ref10}. Generally speaking, due to the lack of a theory for identify the global minimizer(s), many polynomial optimization problems are considered to be \textit{NP}-hard. \\
\indent The canonical duality theory was originated in the late 1980s by Gao and Strang, and has developed significantly in recent years, both theoretically and practically \cite{Ref11,Ref12}. Actually, the canonical duality theory has been successfully applied to solve some special polynomial optimization problems. In \cite{Ref13}, a special polynomial minimization problem called canonical polynomial was completely solved by the canonical duality theory. In \cite{Ref14}, the theory was used to solve a special 8\textit{th} order polynomial minimization problem. Recently, canonical dual solutions to sum of fourth-order polynomials minimization problems have also achieved \cite{Ref15}. This paper aims to solve some general polynomial benchmark problems by using a generalized canonical duality theory. Experimental results show that these polynomial benchmark problems can be solved completely by the canonical duality theory.
\section{A brief review of the canonical duality theory}
\label{sec:1}
Let's consider the following general polynomial optimization problem (primal problem)
\begin{equation}
(\mathcal{P}) : \min_{\mathbf{x} \in \mathbb{R}^n} \Big\{P(\mathbf{x}) = \frac{1}{2} \mathbf{x}^T A \mathbf{x} - \mathbf{x}^T \mathbf{f} + W(\mathbf{x})\Big\},
\end{equation}
where, $A \in \mathbb{R}^{n \times n}$ is a given symmetrical indefinite matrix, $\mathbf{f} \in \mathbb{R}^n$ is a given vector, $W(\mathbf{x}):\mathbb{R}^n \rightarrow \mathbb{R} $ is a general nonconvex $C^2$ function.\\
\indent The main procedures of general methodology of the canonical duality theory can be summarized as the following three steps:\\
\indent Step 1: \textit{Canonical dual transformation}\\
\indent Introducing a nonlinear operator (a G\^{a}teaux differentiable geometrical measure)
\begin{eqnarray}
\bm\xi = \mathrm{\Lambda} (\mathbf{x}): \mathbb{R}^n \rightarrow \mathcal{E}_a \subset \mathbb{R}^m
\end{eqnarray}
and a convex function $V:\mathcal{E}_a \rightarrow \mathbb{R}$ such that $W(\mathbf{x})$ can be recast by $W(\mathbf{x}) = V(\mathrm{\Lambda} (\mathbf{x}))$. Then the primal problem can be rewritten as the canonical form:
\begin{eqnarray}
\min_{\mathbf{x} \in \mathbb{R}^n} \Big\{  P(\mathbf{x}) = V(\mathrm{\Lambda}(\mathbf{x})) - U(\mathbf{x}) \Big\},
\end{eqnarray}
where $U(\mathbf{x}) = - \frac{1}{2}\mathbf{x}^T A \mathbf{x} + \mathbf{x}^T\mathbf{f}$.\\
\indent Step 2: \textit{Generalized complementary function}\\
\indent The dual variable $\bm\varsigma$ to $\bm\xi$ is defined by the duality mapping
\begin{eqnarray}
\bm\varsigma = \nabla V(\bm\xi): \mathcal{E}_a \rightarrow \mathcal{E}^{*}_a \subset \mathbb{R}^m,
\end{eqnarray}
which should be invertible, due to the convexity of $V(\bm\xi)$. Then the Legendre conjugate $V^{\ast}(\bm\varsigma)$ of $V(\bm\xi)$ can be uniquely defined by the Legendre transformation
\begin{eqnarray}
V^{\ast}(\bm\varsigma) = \mathrm{sta} \{ \bm\xi^T\bm\varsigma - V(\bm\xi)|\bm\xi \in \mathcal{E}_a\}
\end{eqnarray}
and the following canonical duality relations hold on $\mathcal{E}_a \times \mathcal{E}^{*}_a$:
\begin{eqnarray}
\bm\varsigma = \nabla V(\bm\xi) \Leftrightarrow \bm\xi  = \nabla V^{\ast}(\bm\varsigma) \Leftrightarrow V(\bm\xi) + V^{\ast}(\bm\varsigma) = \bm\xi^T \bm\varsigma.
\end{eqnarray}
\indent Replacing $W(\mathbf{x}) = V(\mathrm{\Lambda} (\mathbf{x}))$ by $\mathrm{\mathrm{\Lambda }}(\mathbf{x})^T \bm\varsigma - V^{\ast}(\bm\varsigma)$, we obtain the following generalized complementary function:
\begin{eqnarray}
\Xi(\mathbf{x}, \bm\varsigma) = \mathrm{\Lambda}(\mathbf{x})^T\bm\varsigma - V^{\ast}(\bm\varsigma) - U(\mathbf{x}): \mathbb{R}^n \times \mathcal{E}^{*}_a
\rightarrow \mathbb{R}.
\end{eqnarray}
\indent Step 3: \textit{Canonical dual function}\\
\indent By using the generalized complementary function, the canonical dual function $P^d(\bm{\varsigma})$ can be formulated as
\begin{eqnarray}
P^d(\bm{\varsigma}) = \mathrm{sta} \{\Xi(\mathbf{x},\bm{\varsigma})| \mathbf{x} \in \mathbb{R}^n\} = U^{\mathrm{\Lambda}}(\bm\varsigma) - V^{\ast}(\bm\varsigma),
\end{eqnarray}
where $U^{\mathrm{\Lambda}}(\bm\varsigma)$ is defined by
\begin{eqnarray}
U^{\mathrm{\Lambda}}(\bm\varsigma) = \mathrm{sta}\{\mathrm{\Lambda}(\mathbf{x})^T \bm\varsigma - U(\mathbf{x})|\mathbf{x} \in \mathbb{R}^n\}.
\end{eqnarray}
\indent Let $\mathcal{S}_a \subset \mathcal{E}^{*}_a$ be a dual feasible space such that $U^{\mathrm{\Lambda}}(\bm\varsigma)$ is well-defined, and the canonical dual problem can be obtained as
\begin{eqnarray}
(\mathcal{P}^d): \mathrm{sta}\{P^d(\bm{\varsigma})|\bm\varsigma \in \mathcal{S}_a\}.
\end{eqnarray}
\indent \textbf{Theorem 1} (\textbf{Complementary-Dual Principle})\cite{Ref12}.
The problem $(\mathcal{P}^d)$ is canonically dual to the primal problem $(\mathcal{P})$ in the sense that if $(\bar {\mathbf{x}},\bar{\bm{\varsigma}})$ is a critical
point of $\Xi(\mathbf{x}, \bm\varsigma)$, then $\bar{\mathbf{x}}$ is a feasible solution of $(\mathcal{P})$, $\bar{\bm{\varsigma}}$ is a feasible solution of $(\mathcal{P}^d)$, and
\begin{eqnarray}
P(\bar{\mathbf{x}}) = \Xi(\bar{\mathbf{x}}, \bar{\bm{\varsigma}}) = P^d(\bar{\bm{\varsigma}}).
\end{eqnarray}
\indent In many applications, the geometrical operator $\mathrm{\Lambda}(\mathbf{x})$ is usually quadratic
\begin{eqnarray}
\mathrm{\Lambda}(\mathbf{x}) = \{\frac{1}{2}\mathbf{x}^T C_k \mathbf{x} + \mathbf{x}^T \mathbf{b}_k\}: \mathbb{R}^n \rightarrow \mathcal{E}_a \subset \mathbb{R}^m,
\end{eqnarray}
where $C_k \in \mathbb{R}^{n \times n}$ and $\mathbf{b}_k \in \mathbb{R}^n$ are given. In this case, the canonical dual function can be
formulated in the form of 
\begin{eqnarray}
P^d(\bm{\varsigma}) = -\frac{1}{2}\mathbf{F}^T(\bm{\varsigma}) G^{-1}(\bm{\varsigma})\mathbf{F}(\bm{\varsigma}) - V^{*}(\bm{\varsigma}),
\end{eqnarray}
which is well defined on
\begin{eqnarray}
\mathcal{S}_a  = \{\bm{\varsigma} \in \mathbb{R}^m | \mathbf{F}(\bm{\varsigma}) \in \mathcal{C}_{ol}(G(\bm{\varsigma}))\},
\end{eqnarray}
where $G(\bm{\varsigma}) = A + \sum_{k=1}^{m}\varsigma_k C_k$, $\mathbf{F}(\bm{\varsigma}) = \mathbf{f} - \sum_{k=1}^{m}\varsigma_k \mathbf{b}_k$, and $\mathcal{C}_{ol}(G(\bm{\varsigma}))$ denotes the column space of $G(\bm{\varsigma})$.\\
\indent Let the positive domain
\begin{eqnarray}
\mathcal{S}^{+}_a  = \{\bm{\varsigma} \in \mathcal{S}_a | G(\bm{\varsigma}) \succeq 0\}
\end{eqnarray}
where $G(\bm{\varsigma}) \succeq 0$ indicates that $G(\bm{\varsigma})$ is a positive semi-definite matrix.
\\
\indent \textbf{Theorem 2} (\textbf{Global Optimality condition})\cite{Ref12}. Suppose $\bar{\bm{\varsigma}}$ is a critical point of $P^d$ and $\bar{\mathbf{x}} = G^{-1}(\bar{\bm{\varsigma}})F(\bar{\bm{\varsigma}})$. If $\bar{\bm{\varsigma}} \in S^{+}_a$, then $\bar{\bm{\varsigma}}$ is a global maximizer of $(\mathcal{P}^d)$ on $S^{+}_a$ if and only if $\bar{\mathbf{x}}$ is a global minimizer of $(\mathcal{P})$ on $\mathbb{R}^n$, i.e.,
\begin{eqnarray}
P(\mathbf{\bar{x}}) = \min_{\mathbf{x} \in \mathbb{R}^n}P(\mathbf{x}) \Leftrightarrow \max_{\mathbf{\bm{\varsigma}} \in \mathcal{S}^+_a}P^d(\bm{\varsigma}) = P^d(\mathbf{\bar{\bm{\varsigma}}}).
\end{eqnarray}
\indent Some polynomial problems have already been given to testify the effectiveness of canonical duality theory,
but most of them are no more than fourth degree (see \cite{Ref15}). Although some larger degree polynomial problems are completely solved by the same theory, they belong to a special case (see \cite{Ref13,Ref14}). In the next two sections, we are to solve some general polynomial benchmark optimization problems (they are also not obvious to find the global optimum by observation), aiming to expand the use of canonical duality methodology.
\section{Application for Goldstein-Price problem} 
\begin{figure}[htbp!]
{\includegraphics[width=5cm,height=5cm]{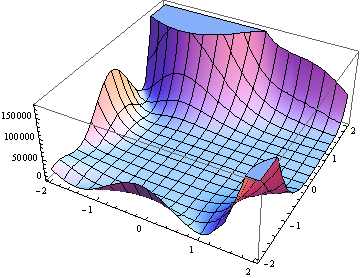}
\includegraphics[width=5cm,height=5cm]{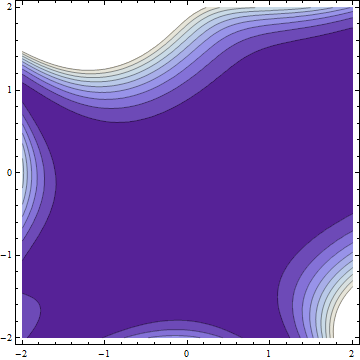}}
\caption{Graph and contour of Goldstein-Price function}
\label{fig1}
\end{figure}

\begin{figure}[htbp!]
{\includegraphics[width=5cm,height=5cm]{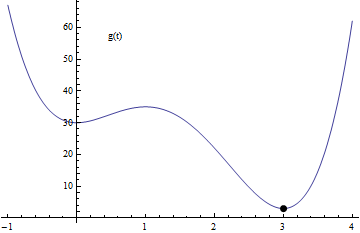}
\includegraphics[width=5cm,height=5cm]{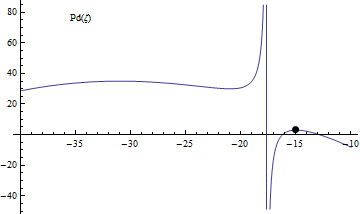}}
\caption{The primal and dual of g(t)}
\label{fig2}
\end{figure}

The Goldstein-Price problem is given in the form of \cite{Ref16}:
\begin{align*}
\min_{x,y} f_1(x,y) &=[1+(x+y+1)^2(19-14x+3x^2-14y+6xy+3y^2)]\times \nonumber \\
&[30+(2x-3y)^2(18-32x+12x^2+48y-36xy+27y^2)]. \nonumber
\end{align*}
\indent The landscape and contour of Goldstein-Price function are given in Fig.\ref{fig1}, and we can find that
there exist a few extrema. Due to the nonconvexity of the problem, it is not easy to find the global minimum.
\\
\indent By using the following linear transformation
\begin{eqnarray}{\label{eqlt}}
\left(
  \begin{array}{c}
    s \\
    t \\
  \end{array}
\right)
= T(x,y)
=
\left(
  \begin{array}{c}
    x +y \\
    2x-3y \\
  \end{array}
\right)
\end{eqnarray}
the Goldstein-Price function can be rewritten to
\begin{eqnarray}
f_1(x,y) = h(s)g(t)
\end{eqnarray}
where
\begin{eqnarray}
h(s) = 1+(s+1)^2(19 - 14s + 3s^2)
\end{eqnarray}
and
\begin{eqnarray}
g(t) = 30 + t^2(18-16t+3t^2)
\end{eqnarray}
\indent \textbf{Proposition 1} Under the linear transformation $T(x,y)$, the Goldstein-Price problem is equivalent to the decoupled problems as follows
\begin{eqnarray}
\min_{x,y}f_1(x,y) = \min_s h(s) \min_t g(t)
\end{eqnarray}
\textit{Proof}. Since the linear transformation in (11) is independent, and $h(s), g(t)$ are bounded below, it is easy to see the proposition follows. ~~~~~~~~~~~~~~~~~~~~$\square$\\
\indent Next, we will solve $\min\limits_s h(s)$ and $\min\limits_{t} g(t)$ separately.\\
\indent For $h(s)$, we can find that $\forall s, h(s) = (1+(s+1)^2((\sqrt{3}s - \frac{7}{\sqrt{3}})^2 + \frac{8}{3})) > 0$, and there is only one critical point $s = -1$ for $h(s)$,
that is to say, $s^{*} = -1$ is the global minimum of $h(s)$.\\
\indent For $g(t)$, we rewrite it to the following canonical form
\begin{eqnarray}
g(t) = V(\mathrm{\Lambda}(t)) - U(t)
\end{eqnarray}
where, $V(\mathrm{\Lambda}(t)) = 3(t^2-\frac{8}{3}t-2)^2 -9(t^2-\frac{8}{3}t-2)$, and $U(t)= - \frac{53}{3}t^2 + 56t$.\\
\indent Introducing a nonlinear operator
\begin{eqnarray}
\xi = \mathrm{\Lambda} (t) = t^2-\frac{8}{3}t-2  = (t- \frac{4}{3})^2 - \frac{34}{9} \geq -\frac{34}{9}
\end{eqnarray}
then
\begin{eqnarray}
V(\xi) = 3\xi^2 - 9\xi,\; \varsigma = 6 \xi - 9 \geq -\frac{95}{3}, \; V^{\ast}(\varsigma) = \frac{(9 + \varsigma)^2}{12},\;
\end{eqnarray}
therefore, we get the generalized complementary function
\begin{eqnarray}
g(\varsigma,t)&= &\Lambda(t)\varsigma - V^{\ast}(\varsigma) - U(t) \nonumber \\
              &= &(t^2-\frac{8}{3}t-2)\varsigma - \frac{\varsigma^2 + 18\varsigma + 81}{12} + \frac{53}{3}t^2 -56t \nonumber \\
              &= &(\varsigma + \frac{53}{3})t^2 - (\frac{8}{3}\varsigma + 56)t - \frac{\varsigma^2 + 18\varsigma +81}{12} - 2\varsigma,
\end{eqnarray}
\indent For a given $\varsigma$, the criticality condition $g_t(\varsigma,t) = 0$ leads to the canonical equilibrium equation
\begin{eqnarray}{\label{eqcee}}
2(\varsigma + \frac{53}{3})t = \frac{8}{3}\varsigma + 56.
\end{eqnarray}
\indent Substituting $t = (\frac{8}{3}\varsigma + 56)/(2(\varsigma + \frac{53}{3}))$ into $g(\varsigma,t)$, we obtain the dual function
\begin{eqnarray}
P^d(\varsigma) = \frac{1}{12} \left(-\varsigma^2-18 \varsigma-81\right)-\frac{\left(\frac{8 \varsigma}{3}+56\right)^2}{4 \left(\varsigma+\frac{53}{3}\right)}-2 \varsigma,
\end{eqnarray}
which is concave in the positive domain
\begin{eqnarray}
S^{+}_a = \{\varsigma | \varsigma + \frac{53}{3} > 0\}.
\end{eqnarray}
The plots of $g(t)$ and $P^d(\varsigma)$ are illustrated in Fig.\ref{fig2}.
\\
\indent Using the sequential quadratic programming method from the Optimization Toolbox within the MATLAB environment for $P^d(\varsigma)$ over $S^{+}_a$, we can get $\varsigma = -15$.  Then we get the corresponding $t^{*} = (\frac{8}{3}\varsigma + 56)/2(\varsigma + \frac{53}{3}) = 3$ by using the canonical equilibrium equation (\ref{eqcee}).\\
\indent By the inverse linear transformation of (\ref{eqlt}), we can obtain the global minimum to $f_1(x,y)$
\begin{eqnarray}
\left(
  \begin{array}{c}
    x^{*} \\
    y^{*} \\
  \end{array}
\right)
= T^{-1}(s^{*},t^{*})
=
\left(
  \begin{array}{c}
    (3s^{*} + t)/5 \\
    (2s^{*} - t)/5 \\
  \end{array}
\right)
=
\left(
  \begin{array}{c}
    0 \\
    -1 \\
  \end{array}
\right), \nonumber
\end{eqnarray}
which is indeed the global minimum with the result given in \cite{Ref16}.
\section{Application for Three Hump Camel Back problem}
\begin{figure}[!htbp]
{\includegraphics[width=5cm,height=5cm]{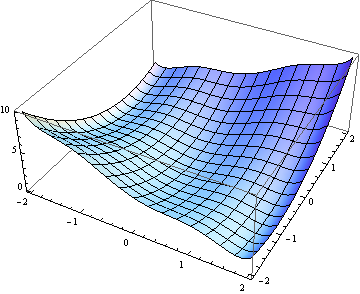}}
{\includegraphics[width=5cm,height=5cm]{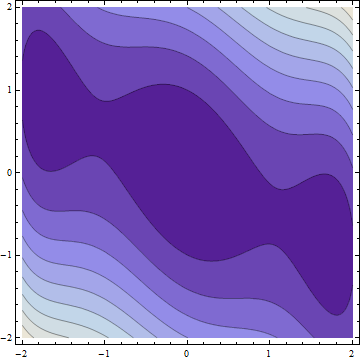}}
\caption{Graph and contour of Three Hump Camel Back function}
\label{fig3}
\end{figure}

The Three Hump Camel Back problem is given in the form of \cite{Ref17}:
\begin{equation*}
\min_{x,y} f_2(x,y)= 2x^2 - 1.05x^4 + \frac{x^6}{6} + xy + y^2.
\end{equation*}
\indent The landscape and contour of Three hump camel back function are given in Fig.\ref{fig3}, and we can find that
there also exist a few extrema. The nonconvexity also makes it difficult to find the global minimum.
\\
\indent Firstly, we rewrite $f_2(x,y)$ to the following canonical form
\begin{eqnarray}
f_2(x,y) = \frac{1}{6}(V_1(\mathrm{\Lambda}_1(x,y)) - U_1(x,y)),
\end{eqnarray}
where, $V_1(\mathrm{\Lambda}_1(x,y)) = (x^3-3.2x)^2$, $U_1(x,y) = - 0.1x^4 - 1.76x^2-6xy-6y^2$.\\
\indent Introducing a nonlinear operator
\begin{eqnarray}
\xi_1 = \mathrm{\Lambda}_1(x) = x^3 - 3.2x,
\end{eqnarray}
then
\begin{eqnarray}
V_1(\xi_1) = \xi^2_1, \; \varsigma_1 = 2\xi_1, \;  V_1^{\ast}(\varsigma_1) = \frac{1}{4}\varsigma^2_1,
\end{eqnarray}
thus, we can get the first generalized complementary function
\begin{eqnarray}
f_2(\varsigma_1,x,y) & = & \frac{1}{6}(\mathrm{\Lambda}_1(x)\varsigma_1 - V_1^{\ast}(\varsigma_1) - U_1(x,y)) \nonumber \\
                     & = & \frac{1}{6}\Big((x^3-3.2x)\varsigma_1 - \frac{1}{4}\varsigma^2_1+0.1x^4 + 1.76x^2 + 6xy + 6y^2 \Big).
\end{eqnarray}
\indent Again, the $f_2(\varsigma_1,x,y)$ can be rewritten to
\begin{eqnarray}
f_2(\varsigma_1,x,y) = \frac{1}{60}(V_2(\mathrm{\Lambda}_2(x,y,\varsigma_1)) - U_2(x,y,\varsigma_1)),
\end{eqnarray}
where, $V_2(\mathrm{\Lambda}_2(x,y,\varsigma_1)) = (x^2+5\varsigma_1 x)^2$, $U_2(x,y,\varsigma_1) = (25\varsigma_1^2 - 17.6)x^2 + 32\varsigma_1x - 60xy - 60y^2 + 2.5\varsigma_1^2$.
\\
\indent Then, we introduce another nonlinear operator
\begin{eqnarray}
\xi_2 = \mathrm{\Lambda}_2(x,\varsigma_1) = x^2 + 5\varsigma_1x,
\end{eqnarray}
thus
\begin{eqnarray}
V_2(\xi_2) = \xi^2_2,\; \varsigma_2 = 2\xi_2, \; V_2^{\ast}(\varsigma_2) =  \frac{1}{4}\varsigma^2_2,
\end{eqnarray}
consequently, we obtain the final generalized complementary function
\begin{eqnarray}
&&f_2(\varsigma_1,\varsigma_2,x,y) = \frac{1}{60}(\mathrm{\Lambda}_2(x,\varsigma_1)\varsigma_2 - V_2^{\ast}(\varsigma_2) - U_2(x,y,\varsigma_1)) \nonumber \\
&&= (\frac{22}{75} - \frac{5}{12}\varsigma^2_1 + \frac{\varsigma_2}{60})x^2 + y^2 +xy + (\frac{1}{12}\varsigma_1\varsigma_2 - \frac{8}{15}\varsigma_1)x -\frac{1}{24}\varsigma^2_1 - \frac{1}{240}\varsigma^2_2~~~~~~
\end{eqnarray}
\indent For given $\varsigma_1$ and $\varsigma_2$, the criticality condition $\nabla\Xi_1(\varsigma_1,\varsigma_2,x,y) = 0$ leads to the following canonical equilibrium equations 
\begin{equation}{\label{eqcees}}
\left\{ \begin{aligned}
         &2(\frac{22}{75} - \frac{5}{12}\varsigma^2_1 + \frac{\varsigma_2}{60})x + y + (\frac{1}{12}\varsigma_1\varsigma_2 - \frac{8}{15}\varsigma_1) = 0 \\
         &2y + x = 0
        \end{aligned}, \right.
\end{equation}
and finally we obtain the canonical dual function
\begin{equation}
P^d(\bm\varsigma) = \frac{-1250 \varsigma_1^4-50 \varsigma_1^2 (31 \varsigma_2-105)+\varsigma_2^2 (5 \varsigma_2+13)}{240 \left(125 \varsigma_1^2-5 \varsigma_2-13\right)},
\end{equation}
which is concave in the positive domain
\begin{eqnarray}
S^{+}_a= \Bigg\{\bm\varsigma\Bigg|
\begin{pmatrix}
\frac{22}{75} - \frac{5}{12}\varsigma^2_1 + \frac{\varsigma_2}{60} & 0.5 \\
0.5 & 1\\
\end{pmatrix}\succeq 0\Bigg\}.
\end{eqnarray}
\indent Using the sequential quadratic programming method from the Optimization Toolbox within the MATLAB environment for $P^d(\varsigma)$ over $S^{+}_a$, we can get $\varsigma_1 = 0, \varsigma_2 = 0$. According to
the canonical equilibrium equations (\ref{eqcees}), we can obtain the corresponding
\begin{eqnarray}
\left(
  \begin{array}{c}
    x^{*} \\
    y^{*} \\
  \end{array}
\right)
=
\left(
  \begin{array}{cc}
    2(\frac{22}{75} - \frac{5}{12}\varsigma^2_1 + \frac{\varsigma_2}{60}) & 1 \\
    1 & 2 \\
  \end{array}
\right)^{-1}
\left(
  \begin{array}{c}
    \frac{8}{15}\varsigma_1 -\frac{1}{12}\varsigma_1\varsigma_2\\
    0 \\
  \end{array}
\right)
=
\left(
  \begin{array}{c}
    0 \\
    0 \\
  \end{array}
\right),
\nonumber
\end{eqnarray}
which is indeed the global minimum with the result given in \cite{Ref17}.
\section{Conclusion} When the canonical duality methodology is applied to a special class of polynomial optimization problem,
it has the ability to solve the class of problem completely. On the other hand, for a general polynomial problem, we can also design appropriate canonical dual transformation to achieve the goal. As for Goldstein-Price problem, we transform it to decoupled minimization problems, and then solve them separately. While for Three hump camel back problem, we can utilize two-level canonical dual transformations. The completely solutions of the general polynomial benchmark functions have witnessed the powerfulness of
the canonical duality methodology again.



\end{document}